\documentclass[12pt,leqno]{article}
\usepackage{amsmath,amssymb,amsthm,amscd}
\numberwithin{equation}{section} \allowdisplaybreaks
\begin{document}
\newtheorem{theorem}{Theorem}[section]
\newtheorem{defin}{Definition}[section]
\newtheorem{prop}{Proposition}[section]
\newtheorem{corol}{Corollary}[section]
\newtheorem{lemma}{Lemma}[section]
\newtheorem{rem}{Remark}[section]
\newtheorem{example}{Example}[section]
\title{Foliation-coupling Dirac structures}
\author{{\small by}\vspace{2mm}\\Izu Vaisman}
\date{}
\maketitle
{\def\thefootnote{*}\footnotetext[1]%
{{\it 2000 Mathematics Subject Classification: 53D17} .
}}
\begin{center} \begin{minipage}{12cm}
A{\footnotesize BSTRACT. We extend the notion of ``coupling with a
foliation" from Poisson to Dirac structures and get the
corresponding generalization of the Vorobiev characterization of
coupling Poisson structures \cite{{Vor},{V04'}}. We show that any
Dirac structure is coupling with the fibers of a tubular
neighborhood of an embedded presymplectic leaf, give new proofs of
the results of Dufour and Wade \cite{DW} on the transversal Poisson
structure, and compute the Vorobiev structure of the total space of
a normal bundle of the leaf. Finally, we use the coupling condition
along a submanifold, instead of a foliation, in order to discuss
submanifolds of a Dirac manifold which have differentiable, induced
Dirac structures. In particular, we get an invariant that reminds
the second fundamental form of a submanifold of a Riemannian
manifold.}
\end{minipage}
\end{center}
\vspace{5mm}
\section{Introduction}
In this paper, the functions, manifolds, bundles, etc. are assumed
to be differentiable of class $C^\infty$. The Dirac structures were
first defined by T. Courant and A. Weinstein \cite{CW} and studied
in Courant's thesis \cite{C}. Dirac structures are important because
they provide a unified view of Poisson and presymplectic structures,
and generalize both. Later, I. Dorfman \cite{Dorf} extended the
notion of a Dirac structure to complexes over Lie algebras. On the
other hand, the bracket used by Courant was extended by Z.-J. Liu,
A. Weinstein and P. Xu \cite{Liu} to a notion of Courant algebroid
and the corresponding generalization of Dirac structures were
introduced and used in \cite{LWX}. An extension of the original
Courant bracket, which is not a Courant algebroid bracket but
includes the Jacobi structures in the scheme, was defined by A\"\i
sa Wade \cite{Wade}.

For the reader's convenience, we recall the general definitions of
\cite{{Liu},{LWX}} in a slightly different form. A {\it Courant
algebroid} is a vector bundle $p:C\rightarrow M$ endowed with a
non-degenerate, pseudo-Euclidean metric $g\in\Gamma\odot^2C$
($\Gamma$ denotes spaces of cross sections and $\odot$ denotes
symmetric tensor product), a morphism $\rho:C\rightarrow TM$ (the
{\it anchor}) and a skew-symmetric bracket $[\,,\,]_C:\Gamma
C\times\Gamma C\rightarrow\Gamma C$ such that:
$$\begin{array}{rl}i)&\hspace*{3mm}\rho[c_1,c_2]_C=[\rho c_1,\rho
c_2]_{TM},\vspace{2mm}\\
ii)&\hspace*{3mm}\sum_{Cycl(1,2,3)}[[c_1,c_2]_C,c_3]_C =
\frac{1}{3}\partial\{\sum_{Cycl(1,2,3)}g([c_1,c_2]_C,c_3)\},\vspace{2mm}\\
iii)&\hspace*{3mm}(\rho c)\{g(c_1,c_2)\} = g([c,c_1]_C+\partial
g(c,c_1),c_2) + g(c_1,[c,c_2]_C+\partial g(c,c_2)),\end{array}$$
where $c_a\in\Gamma C$ $(a=1,2,3)$ and, if $f\in C^\infty(M)$,
$\partial f=(1/2)\sharp_g\circ^t\hspace{-2pt}\rho(df)$ ($^t$ denotes
transposition and the ``musical morphisms" are defined like in
Riemannian geometry), equivalently,
\begin{equation}\label{exprpartial} g(c,\partial f)=\frac{1}{2}
(\rho c)f.\end{equation} (Further basic properties of Courant
algebroids may be found, for example, in \cite{{Royt},{V04}}.)

The most important Courant algebroids are the so called doubles of
Lie bialgebroids \cite{{Liu},{LWX}}. We describe them by means of
the notion of a {\it para-Hermitian structure} on a pseudo-Euclidean
bundle $(C,g)$ (e.g., \cite{CFG}) that is, a bundle morphism
$F:C\rightarrow C$ such that
\begin{equation}\label{parastr}F^2=Id.,\;
g(Fc_1,Fc_2)=-g(c_1,c_2)\hspace{5mm}(\forall c_{1,2}\in\Gamma
C).\end{equation} A para-Hermitian vector bundle $(C,g,F)$
decomposes as $C=C_{+}\oplus C_{-}$, where the components are the
$(\pm1)$-eigenspaces of $F$, respectively, and the projections onto
these components are
\begin{equation}\label{parapm}F_{\pm}=\frac{1}{2}(Id.\pm F).\end{equation}
Moreover, $C_{\pm}$ are maximal isotropic with respect to $g$, the
latter must be {\it neutral} (i.e., of signature zero), and one has
isomorphisms $\flat_g:C_{\pm}\rightarrow C_{\mp}^*$ (the star
denotes the dual bundle). Furthermore, the bundle also has the
non-degenerate $2$-form
\begin{equation}\label{omega}
\omega(c_1,c_2)=
g(c_1,Fc_2)\hspace{1cm}(\omega(Fc_1,Fc_2)=-\omega(c_1,c_2))
\end{equation}
and the subbundles $C_{\pm}$ are $\omega$-Lagrangian. All these
facts apply to para-Hermitian vector spaces, which is the case where
the basis $M$ is a point.

A Courant algebroid $(C,g,\rho,[\,,\,]_C)$ will be called {\it
para-Hermitian} if it is endowed with a para-Hermitian structure $F$
such that the subbundles $C_{\pm}$ are closed with respect to the
bracket $[\,,\,]_C$ i.e., $\forall c_1,c_2\in\Gamma C$, one has
\begin{equation}\label{inchidere}
F_{-}[F_{+}c_1,F_{+}c_2]_C=0,\;F_{+}[F_{-}c_1,F_{-}c_2]_C=0.
\end{equation} Taking into account the expression of $F_{\pm}$ and the
property $F^2=Id.$, we see that the two conditions above are
equivalent with the following single condition
\begin{equation}\label{integrabil}
[Fc_1,Fc_2]_C-F[Fc_1,c_2]_C-F[c_1,Fc_2]_C+[c_1,c_2]_C=0,
\end{equation} which will
be called the {\it integrability condition} of $F$, because this is
the integrability condition of a para-Hermitian structure on the
tangent bundle of a manifold.

Since the subbundles $C_{\pm}$ are $g$-isotropic, the Courant
algebroid axioms imply that the vector bundle structure of a
para-Hermitian Courant algebroid is that of a direct sum of two dual
Lie algebroids of anchors $\rho\circ F_{+},\rho\circ F_{-}$.
Moreover, the Lie algebroid brackets of $C_{\pm}$ together with $F$
and $g$ determine the Courant bracket of $C$. Indeed, from
(\ref{parapm}) and (\ref{inchidere}), it follows that
\begin{equation}\label{crosetcuF} \begin{array}{l}
F_{+}[c_1,c_2]_C= [F_{+}c_1,F_{+}c_2]_C+ F_{+}([F_{+}c_1,F_{-}c_2]_C
+ [F_{-}c_1,F_{+}c_2]_C),\vspace{2mm}\\ F_{-}[c_1,c_2]_C=
[F_{-}c_1,F_{-}c_2]_C+ F_{-}([F_{+}c_1,F_{-}c_2]_C +
[F_{-}c_1,F_{+}c_2]_C).
\end{array}\end{equation} Then, by writing down axiom iii)
of the definition of a Courant algebroid for triples
$(F_{+}c,F_{+}c_1,F_{-}c_2)$, $(F_{-}c_2,F_{-}c_1,F_{+}c)$, instead
of $(c,c_1,c_2)$, using (\ref{exprpartial}) and conveniently
permuting $(c,c_1,c_2)$, we get the formulas
\begin{equation}\label{crosetuldepm} \begin{array}{c}
g(F_{+}c,[F_{+}c_1,F_{-}c_2]_C) = g(F_{-}c_2,[F_{+}c,F_{+}c_1]_C)
\vspace{2mm}\\

+ (\rho F_{+}c_1)(g(F_{+}c,F_{-}c_2)) -\frac{1}{2}(\rho
F_{+}c)(g(F_{+}c_1,F_{-}c_2)),\vspace{2mm}\\

g(F_{-}c,[F_{+}c_1,F_{-}c_2]_C) = -g(F_{+}c_1,[F_{-}c,F_{-}c_2]_C)
\vspace{2mm}\\

+ (\rho F_{-}c_2)(g(F_{-}c,F_{+}c_1)) +\frac{1}{2}(\rho
F_{-}c)(g(F_{+}c_1,F_{-}c_2)).
\end{array}\end{equation} These formulas define the brackets
$[F_{+}c_1,F_{-}c_2]_C$, which, together with (\ref{crosetcuF}),
proves the previous assertion.

Accordingly, one can see that the notion of a para-Hermitian Courant
algebroid is the same as that of the double of a Lie bialgebroid
\cite{Liu}.

An {\it almost Dirac structure} of the para-Hermitian Courant
algebroid $C$ is a maximal $g$-isotropic subbundle $L$ of $C$
\cite{LWX}. The isotropy property may be expressed by
\begin{equation}\label{isotropialuiL}
g(F_{+}l_1,F_{-}l_2)+ g(F_{-}l_1,F_{+}l_2)=0,\hspace{5mm} \forall
l_1,l_2\in\Gamma L.\end{equation}

The algebraic properties of almost Dirac structures were discussed
in \cite{C}. Put $p_{\pm}=F_{\pm}|_L$. Then,
$ker\,p_{\pm}=C_{\mp}\cap L$, and we get the subbundles
$$L_{\pm}=im\,p_{\pm}\approx L/{C_{\mp}\cap
L}.$$ It is easy to see that $ker(\omega|_L)=(C_{+}\cap L)\oplus
(C_{-}\cap L)$, where $\omega$ is the $2$-form (\ref{omega}). Hence,
the subbundles $L_{\pm}$ have induced $2$-forms $\omega_{\pm}^L$
defined by
\begin{equation}\label{eqptomegapm}
\omega_{+}^L(F_{+}l_1,F_{+}l_2)= \omega_{-}^L(F_{-}l_1,F_{-}l_2)
=\omega(l_1,l_2)\end{equation}
$$=2g(F_{-}l_1,F_{+}l_2)=-2g(F_{+}l_1,F_{-}l_2),\hspace{5mm}
l_{1},l_{2}\in\Gamma L.$$

It is possible to reconstruct $L$ from each of the pairs
$(L_{\pm},\omega_{\pm}^L)$. Namely, with (\ref{eqptomegapm}) one
gets
\begin{equation}\label{reconstrL+}
\begin{array}{r} L=\{c\;/\;F_{+}(c) \in
L_{+},\;g(F_{-}(c),u)=\frac{1}{2}\omega_{+}^L(F_{+}(c),u),\forall
u\in L_{+}\},\vspace{2mm}\\ L=\{c\;/\;F_{-}(c) \in
L_{-},\;g(F_{+}(c),v)=-\frac{1}{2}\omega_{-}^L(F_{-}(c),v),\forall
v\in L_{-}\}. \end{array} \end{equation}

In particular, if $L_{+}=C_{+}$, $L$ is determined by the $2$-form
$\omega_{+}$ on $C_{+}$ and it may be called an {\it almost
presymplectic structure}. In this case, the first formula
(\ref{reconstrL+}) shows that $L$ may be identified with the graph
of the mapping $(1/2)\sharp_g\circ\flat_\omega:C_{+}\rightarrow
C_{-}$. If $L_{-}=C_{-}$, $L$ is determined by the $2$-form
$\omega_{-}$ on $C_{-}$, it may be called an {\it almost Poisson
structure} and it is the graph of the mapping
$-(1/2)\sharp_g\circ\flat_{\omega_{-}}:C_{-}\rightarrow C_{+}$. The
condition $L_{+}=C_{+}$ is equivalent with the surjectivity of
$p_{+}$, i.e., with $ker\,p_{+}=C_-\cap L=\{0\}$, and this latter
condition also characterizes the almost presymplectic case.
Similarly, the almost Poisson case is also
characterized by $C_{+}\cap L=\{0\}$.\\

Finally, a {\it Dirac structure} is an almost Dirac structure which
is closed with respect to the bracket $[\,,\,]_C$. Equivalently,
$L\subseteq C$ is a Dirac structure if it is maximal isotropic and
$\forall l_a\in\Gamma L$ $(a=1,2,3)$ one has
\begin{equation}\label{condinchidere}
g([l_1,l_2]_C,l_3)=0.\end{equation} From the axioms of the Courant
algebroids it follows that if $L$ is a Dirac structure
then $(L,\rho|_L,[\,,\,]_C)$ is a Lie algebroid.\\

In this paper, we will only be interested in the classical Courant
case \cite{C}. That is $C=TM\oplus T^*M$ with $\rho(X,\alpha)=X$,
\begin{equation}\label{gFinC}
g((X,\alpha),(Y,\beta))=\frac{1}{2}(\beta(X)+\alpha(Y)),\;F(X,\alpha)=(X,-\alpha),
\end{equation} therefore, \begin{equation}\label{omegainC}
\omega((X,\alpha),(Y,\beta))= \frac{1}{2}(\alpha(Y)-\beta(X)),
\end{equation} and with the bracket
\begin{equation}\label{crosetinC} [(X,\alpha),(Y,\beta)] = ([X,Y],
L_X\beta-L_Y\alpha+d(\omega((X,\alpha),(Y,\beta))))\end{equation}
$$= ([X,Y],i(X)d\beta-i(Y)d\alpha +
\frac{1}{2}d(\beta(X)-\alpha(Y))).$$ In the previous formulas, $X,Y$
are vector fields and $\alpha,\beta$ are $1$-forms on the
differentiable manifold $M$, and the bracket of vector fields is the
usual Lie bracket. Notice that $C_{+}=TM$, $C_{-}=T^*M$ and the
Courant bracket reduces to zero on $T^*M$.

Then, a (almost) Dirac structure of $TM\oplus T^*M$ is called a
(almost) Dirac structure on the manifold $M$. By the first formula
(\ref{reconstrL+}), an almost Dirac structure $L$ of $M$ is
determined by a generalized distribution $L_{+}\subseteq TM$ endowed
with a $2$-form $\omega^L_{+}$, namely:
\begin{equation}\label{reconstrincazclasic}
L=\{(X,\alpha)\,/\,X\in L_{+}\,\&\,\alpha|_{L_{+}}=
\flat_{\omega^L_{+}}X\}.\end{equation} By a technical computation,
it follows from (\ref{reconstrincazclasic}) that $L$ is a Dirac
structure iff $L^+$ is integrable and the form $\omega^L_{+}$ is
closed on the leaves of $L^+$ \cite{C}. Accordingly, a Dirac
structure on $M$ is equivalent with a generalized foliation with
presymplectic leaves where the presymplectic form depends
differentiably of the leaves. If the leaves are symplectic we have a
Poisson structure, and if the leaves are the connected components of
$M$ we have a presymplectic structure (of a non-constant rank) on
$M$.

Jacobi structures on a manifold $M$ may be seen as a particular case
of Dirac structures on $M\times\mathbb{R}$. Namely, a Jacobi
structure on $M$ is equivalent with a {\it Poisson homogeneous
structure} on $M\times\mathbb{R}$ (e.g., \cite{DLM}). We recall that
the Poisson structure defined by the bivector field $P$ is
homogeneous if there exists a vector field $Z$ such that $P+L_ZP=0$.
This is equivalent with the fact that the Dirac structure
$L(P)=\{(\sharp_P\theta,\theta)\,/\,\theta\in T^*M\}$ is such that
$\forall(X,\theta)\in L(P)$ one has $(X+[Z,X],L_Z\theta)\in L(P)$.
The latter property may be attributed to a general Dirac structure
$L$, thus, producing the notion of a general {\it homogeneous Dirac
structure}. A more sophisticated way to see Jacobi structures as
Dirac was proposed in \cite{Wade}.

Furthermore, we will be interested in the case where the manifold
$M$ is also endowed with a regular foliation $\mathcal{F}$, and our
aim is to extend the notion of {\it $\mathcal{F}$-coupling} from
Poisson structures to Dirac structures. Poisson structures coupling
with a fibration were studied by Vorobiev \cite{Vor} then, extended
to foliated manifolds in \cite{V04'}. They proved to be important in
the study of the geometry of a Poisson structure in the neighborhood
of an embedded symplectic leaf \cite{Vor}. In \cite{DW}, Dufour and
Wade study a Dirac structure in the neighborhood of a presymplectic
leaf and (in our terms) show that the structure is coupling with
respect to the fibers of a tubular neighborhood. In the present
paper we will define the coupling property of a Dirac structure with
respect to an arbitrary foliation and extend Vorobiev's results. In
particular, we will give geometric proofs of some of the results of
\cite{DW}.

Since this paper is a continuation of \cite{V04'}, and in order to
avoid repetition, we will use the same notation for everything
related with the foliation. In particular, we assume that $dim\,M=n,
dim\,\mathcal{F}=p, q=n-p$, and we denote by
$\Omega^*(M),\mathcal{V}^*(M)$ the spaces of differential forms and
multivector fields on $M$. We will use a normal bundle $H$, i.e.,
$TM=H\oplus F$, $F=T\mathcal{F}$, and $T^*M=H^*\oplus F^*$ for the
dual bundles $H^*=ann\,F,F^*=ann\,H$ ($ann$ denotes the annihilator
space). We will also use the corresponding bigrading of differential
forms and multivector fields and the induced decomposition
\begin{equation}\label{descompd} d=d'_{1,0} + d''_{0,1}
+\partial_{2,-1}\end{equation} of the exterior differential.
\section{Coupling Dirac structures}
Let $(M,\mathcal{F})$ be a foliated manifold as described at the end
of Section 1. From \cite{V04'}, we recall that a bivector field
$P\in\mathcal{V}^2(M)$ is $\mathcal{F}$-almost coupling via the
normal bundle $H$ if $P=P'_{2,0}+P''_{0,2}$, where  the indices
denote the bidegree, i.e., $P'\in\Gamma\wedge^2H,
P''\in\Gamma\wedge^2F$. In this case, $P$ satisfies the Poisson
condition $[P,P]=0$ iff the following four conditions hold:
\begin{equation}\label{Poissonac} \begin{array}{l}
(L_{\sharp_{P'}\gamma}P')(\alpha,\beta)=
d'\gamma(\sharp_{P'}\alpha,\sharp_{P'}\beta),\vspace{2mm}\\
(L_{\sharp_{P'}\lambda}P')(\alpha,\beta)=
-\lambda([\sharp_{P'}\alpha,\sharp_{P'}\beta]),\vspace{2mm}\\
(L_{\sharp_{P'}\gamma}P'')(\lambda,\mu)=0,\vspace{2mm}\\
(L_{\sharp_{P''}\nu}P')(\lambda,\mu)=d''\nu(\sharp_{P''}\lambda,\sharp_{P''}\mu),
\end{array}\end{equation}
where $\alpha,\beta,\gamma\in\Omega^{1,0}(M),\lambda,\mu,\nu
\in\Omega^{0,1}(M)$.

Accordingly, the generalization of the almost coupling condition has
to ask for a decomposition of the Dirac structure into an
$F$-component and an $H$-component.
\begin{defin}\label{precD} {\rm Let $L\subseteq TM\oplus T^*M$ be
a maximal isotropic subbundle. Denote
\begin{equation}\label{componenteLHF} L_H=
L\cap(H\oplus H^*),\; L_F=L\cap(F\oplus F^*).
\end{equation} Then, the almost Dirac structure $L$ is {\it
$\mathcal{F}$-almost coupling via $H$} if}
\begin{equation}\label{defacD} L=L_H\oplus L_F. \end{equation}
\end{defin}

Therefore, $L$ is almost coupling iff $(Z,\theta)\in L$ is
equivalent with $(X,\alpha)$, $(Y,\lambda)\in L$, where
$Z=X+Y,\,\theta=\alpha+\lambda$, $X\in\Gamma H, Y\in\Gamma
F,\alpha\in \Gamma H^*,\lambda\in\Gamma F^*$. Another important
observation that follows from (\ref{defacD}) is that $L_H,L_F$ are
maximal isotropic in $H\oplus H^*,F\oplus F^*$, respectively, for
the metrics induced by $g$ of (\ref{gFinC}). It follows easily that
the bivector field $P$ is $\mathcal{F}$-almost coupling via $H$ iff
the subbundle $L(P)$ satisfies condition (\ref{defacD}). If
$L=L(\tau)$ is an almost presymplectic structure defined by a
$2$-form $\tau$, almost coupling via $H$ holds iff
$\tau=\tau'_{2,0}+\tau''_{0,2}$, where, again, indices denote the
bidegree.

In the almost coupling situation, the integrability condition of a
maximally isotropic subbundle $L\subseteq TM\oplus T^*M$ extends
conditions (\ref{Poissonac}).
\begin{prop}\label{printegrabacD} The $\mathcal{F}$-almost coupling,
almost Dirac structure $L\subseteq TM\oplus T^*M$ is a Dirac
structure iff, $\forall (X,\alpha)\in\Gamma L_H, \forall
(Y,\lambda)\in\Gamma L_F$, the following four conditions hold:
\begin{equation}\label{integrabacD} \begin{array}{l}
\sum_{Cycl(1,2,3)}\{X_1(\alpha_2(X_3))+\alpha_1([X_2,X_3])\} =
0,\vspace{2mm}\\
(L_Y\alpha_2)(X_1)+\alpha_1([Y,X_2])=\lambda([X_1,X_2]),
\vspace{2mm}\\ (L_X\lambda_1)(Y_2)+\lambda_2([X,Y_1])=0,
\vspace{2mm}\\ ([Y_1,Y_2], i(Y_1)d''\lambda_2 - i(Y_2)d''\lambda_1 +
d''(\lambda_2(Y_1)))\in L_F. \end{array}\end{equation}\end{prop}
\begin{proof} Since $L$ is isotropic, by (\ref{isotropialuiL}),
$\forall(Z_a,\theta_a)\in L$ $(a=1,2)$ we have
\begin{equation}\label{coroliso}
\theta_1(Z_2)+\theta_2(Z_1)=0.\end{equation} Then, since $F$ is
involutive and using the decomposition (\ref{descompd}), the second
expression (\ref{crosetinC}) of the Courant bracket yields
\begin{equation}\label{tr-tr} [(X_1,\alpha_1),(X_2,\alpha_2)] =
(pr_H[X_1,X_2], i(X_1)d'\alpha_2-i(X_2)d'\alpha_1 \end{equation}
$$+d'(\alpha_2(X_1)))
+ (pr_F[X_1,X_2], i(X_1)d''\alpha_2-i(X_2)d''\alpha_1 +
d''(\alpha_2(X_1))),$$
\begin{equation}\label{tr-tg} [(X,\alpha),(Y,\lambda)] =
(pr_H[X,Y], i(X)\partial\lambda - i(Y)d''\alpha)\end{equation}$$ +
(pr_F[X,Y], i(X)d'\lambda),$$
\begin{equation}\label{tg-tg} [(Y_1,\lambda_1),(Y_2,\lambda_2)]=
(0, i(Y_1)d'\lambda_2 - i(Y_2)d'\lambda_1 + d'(\lambda_2(Y_1)))
\end{equation} $$+ ([Y_1,Y_2], i(Y_1)d''\lambda_2 - i(Y_2)d''\lambda_1 +
d''(\lambda_2(Y_1))), $$ where $pr$ denotes natural projections and
the terms are $H\oplus H^*$ and $F\oplus F^*$ components,
respectively. In the almost coupling situation, integrability means
that these components always belong to $L_H,L_F$, respectively.

The second term of (\ref{tg-tg}) yields the fourth condition
(\ref{integrabacD}). Since $d''$ is exterior differentiation along
the leaves of $\mathcal{F}$, this condition is equivalent with the
fact that $L_F$ consists of Dirac structures on the leaves of
$\mathcal{F}$; we will say that $L_F$ is a {\it leaf-tangent Dirac
structure} on $(M,\mathcal{F})$.

By maximal isotropy, the $(H\oplus H^*)$-component of (\ref{tg-tg})
belongs to $L_H$ iff, $\forall(X,\alpha)\in\Gamma L_H$, the $1$-form
of the first term of the right hand side of (\ref{tg-tg}) vanishes
on $X$. The result of this evaluation exactly is the third condition
(\ref{integrabacD}).

The terms of the decompositions (\ref{tr-tg}), (\ref{tr-tr}) will be
treated in a similar way, i.e., using maximal isotropy and
evaluations of exterior differentials. The computations show that
the condition provided by the $(F\oplus F^*)$-component of
(\ref{tr-tg}) is again the third condition (\ref{integrabacD}), and
the condition provided by the $(H\oplus H^*)$-component of
(\ref{tr-tg}) is the second condition (\ref{integrabacD}). Then, the
condition provided by the $(F\oplus F^*)$-component of (\ref{tr-tr})
is again the second condition (\ref{integrabacD}), and the condition
provided by the $(H\oplus H^*)$-component of (\ref{tr-tr}) is the
first condition (\ref{integrabacD}).
\end{proof}
\begin{rem}\label{obscazPsiDW}
{\rm With a few computations, one can see that the Poisson
conditions (\ref{Poissonac}) for an almost coupling bivector field
are exactly the Dirac conditions (\ref{integrabacD}) for the
subbundle $L(P)$, and in the same order. If only the component $L_F$
is of the almost Poisson type, i.e., the graph of a bivector field
$\Pi\in\Gamma\wedge^2F$, the last formula (\ref{integrabacD}) means
that $\Pi$ must be a leaf-tangent Poisson structure of $\mathcal{F}$
\cite{V04'} and, by putting $Y_a=\sharp_{\Pi}\lambda_a$ $(a=1,2)$ in
the third formula (\ref{integrabacD}), the latter becomes
\begin{equation}\label{cond3inDW}(L_X\Pi)(\lambda_1,\lambda_2)=0.
\end{equation}}\end{rem}\vspace*{1mm}

Now, on a foliated manifold $(M,\mathcal{F})$, a bivector field $P$
is $\mathcal{F}$-coupling if $\sharp_P(ann\,F)$ is a normal bundle
of the foliation $\mathcal{F}$. In order to extend this notion, with
any maximal isotropic subbundle  $L\subseteq TM\oplus T^*M$, we
associate the possibly non-differentiable, generalized distribution
of $M$ defined by
\begin{equation}\label{distribHLF}H_x(L,\mathcal{F}) =
\{Z\in T_xM\,/\,\exists\alpha\in ann\,F_x\& (Z,\alpha)\in
L\}\hspace{3mm}(x\in M).\end{equation} Then, state the following:
\begin{defin} \label{couplingptD} {\rm The almost Dirac structure
$L$ is $\mathcal{F}$-{\it coupling} if the distribution
$H=H(L,\mathcal{F})$ is normal to the foliation $\mathcal{F}$ at
each point $x\in M$.}
\end{defin}
\begin{prop}\label{cimplicaac} If the subbundle $L$ is
$\mathcal{F}$-coupling, $\forall x\in M$, $L_x$ is
$\mathcal{F}$-almost coupling at $x$ via $H_x=H_x(L,\mathcal{F})$.
Furthermore, the $(H\oplus H^*)$-component of $L_x$ is the graph of
a mapping $\flat_{\sigma_x}:H_x\rightarrow H^*_x$ defined by some
$\sigma_x\in\wedge^2(ann\,F_x)$, and the $(F\oplus F^*)$-component
of $L_x$ is the graph of a mapping $\sharp_{\Pi_x}:F^*_x\rightarrow
F_x$ defined by some $\Pi_x\in\wedge^2F_x$. Moreover,
$H=H(L,\mathcal{F})$ is a differentiable, normal bundle of
$\mathcal{F}$ such that $L$ is $\mathcal{F}$-almost coupling via
$H$, and the global cross sections
$\sigma\in\Gamma\wedge^2(ann\,F),\Pi\in\Gamma\wedge^2F$ are
differentiable.
\end{prop}
\begin{proof} The following considerations are at a fixed point
$x\in M$, which we do not include in the notation. With
$H=H(L,\mathcal{F})$ as the normal space of $\mathcal{F}$ at $x$,
take $(Z,\theta)\in L$ and decompose $Z=X+Y$, $X\in H, Y\in F$. By
the definition of $H$, $\exists\alpha\in H^*$ such that
$(X,\alpha)\in L$, and we get a decomposition
\begin{equation}\label{descompaux}
(Z,\lambda)=(X,\alpha)+(Y,\theta-\alpha),\end{equation} where the
terms belong to $L$. Then, $\forall X'\in H$ with a corresponding
covector $\alpha'$ such that $(X',\alpha')\in L$, (\ref{coroliso})
implies
$$(\theta-\alpha)(X')=-\alpha'(Y)=0.$$ Hence $\theta-\alpha\in
F^*$ and (\ref{descompaux}) implies the almost coupling property
(\ref{defacD}) at $x$.

Furthermore, if $(X,\alpha),(X,\alpha')\in L$, where
$\alpha,\alpha'\in H^*$, we get $(0,\alpha'-\alpha)\in L$ and the
isotropy of $L$ together with the coupling hypothesis imply
$\alpha'=\alpha$. Therefore, $\forall X\in H$, the covector
$\alpha\in ann\,F$ such that $(X,\alpha)\in L$ is unique and $L_H$
is the graph of a morphism $\flat_\sigma$, $\sigma\in \wedge^2H^*$.
Notice also that the uniqueness of $\alpha$ is equivalent with
\begin{equation} \label{intersectLann} L\cap ann\,F=\{0\}.
\end{equation} On the other hand, the definition of $H$ implies
$L\cap F\subseteq H$, therefore, in the coupling case, $L\cap
F=\{0\}$ and (see the Introduction) $L_F$ must be of the almost
Poisson type, whence the existence of $\Pi$.

Finally, we will prove the differentiability of the distribution
$H(L,\mathcal{F})$. For this purpose, let us consider the subspaces
\begin{equation}\label{supraH} \tilde H(L,\mathcal{F})=\{(Z,\alpha)\in L\,/\,
\alpha\in ann\,F\}=L\cap[TM\oplus(ann\,F)]\end{equation} at each
point of $M$. Then, $ker\,p_{+}|_{\tilde H(L,\mathcal{F})}=L\cap
ann\,F$ and
\begin{equation}\label{HcusupraH} H(L,\mathcal{F})=p_{+}(\tilde H(L,\mathcal{F}))
\approx \tilde H(L,\mathcal{F})/(L\cap ann\,F).\end{equation} In the
coupling case, because of (\ref{intersectLann}), $p_{+}|_{\tilde
H(L,\mathcal{F})}$ is an isomorphism and we are done if we prove the
differentiability of $\tilde H(L,\mathcal{F})$. In a neighborhood of
a point $x$, let $\{\alpha^a\}$, $\{\lambda^u\}$ and
$\{(Z_i,\theta_i=\theta_{ia}\alpha^a+\theta_{iu}\lambda^u)\}$ be
differentiable, local bases of $H^*$, $F^*$ and $L$, respectively.
Then,
$$\tilde H(L,\mathcal{F})=\{(\xi^iZ_i,\xi^i\theta_i)\,/\,\theta_{iu}\xi^i=0\}$$
(here and in the whole paper we use the Einstein summation
convention), and we see that $\tilde H(L,\mathcal{F})$ is locally
generated by fundamental solutions of a linear, homogeneous system
of equations with differentiable coefficients. But, if the rank of
the latter is constant (and under the coupling hypothesis the rank
is $n-q$), differentiable, fundamental solutions exist. Of course,
the differentiability of $H$ also implies the differentiability of
the $2$-form $\sigma$ and of the bivector field $\Pi$.
\end{proof}

Hereafter, in the coupling situation we will use only
$H(L,\mathcal{F})$ as the normal bundle and shortly denote it by
$H$. The coupling situation is interesting precisely because it
provides a canonical normal bundle of $\mathcal{F}$.
\begin{prop}\label{geomdata1} An $\mathcal{F}$-coupling, almost
Dirac structure $L\subseteq TM\oplus T^*M$ is equivalent with a
triple $(H,\sigma,\Pi)$ where $H$ is a normal bundle of the
foliation $\mathcal{F}$, $\sigma\in\Gamma\wedge^2(ann\,F)$ and
$\Pi\in\Gamma\wedge^2F$.\end{prop}
\begin{proof} We have already derived the triple from the
subbundle $L$. Such a triple is called a system of {\it geometric
data} \cite{{Vor},{DW}}. Conversely, if the geometric data are
given, we reconstruct $L=L_H\oplus L_F$ by defining $L_H$ as the
graph of $\flat_\sigma$ and $L_F$ as the graph of $\sharp_\Pi$. In
other words, we have
\begin{equation}\label{Lreconstruit}
L=\{(X,\flat_\sigma X)+(\sharp_{\Pi}\lambda,\lambda)\,/\,X\in
H,\,\lambda\in F^*\}.\end{equation}\end{proof}
\begin{corol}\label{Fernandes} On $(M,\mathcal{F})$,
the almost Dirac structure $L$ is $\mathcal{F}$-coupling iff
\begin{equation}\label{Rui} L\cap(F\oplus ann\,F)=\{0\}.
\end{equation} \end{corol} \begin{proof} Condition (\ref{Rui}) is
an immediate consequence of (\ref{Lreconstruit}). Conversely,
(\ref{Rui}) implies $H\cap F=\{0\}$ and $H\approx\tilde H$. On the
other hand, by looking at dimensions, (\ref{Rui}) also implies
$L\oplus(F\oplus ann\,F)=TM\oplus T^*M$, therefore, $L+(TM\oplus
ann\,F)=TM\oplus T^*M$, and from (\ref{supraH}) we get
$dim\,H=\dim\,\tilde H=q$. \end{proof}
\begin{rem}\label{obs2formac} {\rm An almost Poisson structure $L(P)$
defined by the bivector field $P$ is coupling iff there exists a
normal bundle $H$ of the foliation $\mathcal{F}$ that yields
$P=P'_{2,0}+P''_{0,2}$ where $P'$ is non degenerate \cite{Vor}. In
the case of an almost presymplectic structure $L(\tau)$ defined by a
$2$-form $\tau$, $H(L(\tau),\mathcal{F})$ is the $\tau$-orthogonal
distribution of $F$ and the coupling condition holds iff the former
is a complementary distribution of the latter. Equivalently,
$L(\tau)$ is $\mathcal{F}$-coupling iff there exists a normal bundle
$H$ that yields a decomposition
\begin{equation}\label{descompsigma}\tau=\tau'_{2,0}+\tau''_{0,2},
\end{equation} where
$\tau''$ is non degenerate.}\end{rem}
\begin{rem}\label{cocoupling} {\rm One can also define the notion
of an $\mathcal{F}$-coupling Dirac structure in a dual way. Namely,
for any almost Dirac structure $L\subseteq TM\oplus T^*M$ of the
foliated manifold $(M,\mathcal{F})$, we may define the {\it
generalized codistribution} (a field of subspaces of the fibers of
$T^*M$ with a varying dimension)
\begin{equation}\label{campconormal} K^*_x=K^*_x(L,\mathcal{F}) =
\{\theta\in T^*_xM\,/\,\exists Y\in F_x\,\&\,(Y,\theta)\in L_x\}
\hspace{3mm}(x\in M).\end{equation} $K^*$ may not be differentiable,
i.e., it may not have local generators defined by differentiable
$1$-forms. Then, it follows that $L$ is $\mathcal{F}$-coupling iff
\begin{equation}\label{eqco-coupling}T^*M=(ann\,F)\oplus K^*.
\end{equation} Indeed, by dualizing the proof of Proposition
\ref{cimplicaac}, we see that condition (\ref{eqco-coupling}) also
obliges $L$ to be of the form (\ref{Lreconstruit}). Notice also
that, in the coupling case, the decomposition (\ref{eqco-coupling})
is the dual of $TM=H\oplus F$ for $H$ given by
(\ref{distribHLF}).}\end{rem}

The integrability conditions of a coupling Dirac structure may also
be expressed by means of the associated geometric data like in the
Poisson case \cite{{Vor}, {V04'}}.
\begin{prop}\label{integrabcugd} An $\mathcal{F}$-coupling almost Dirac
structure $L\subseteq TM\oplus T^*M$ of a foliated manifold
$(M,\mathcal{F})$ is a Dirac structure iff its associated geometric
data $(H,\sigma,\Pi)$ satisfy the following
conditions:\\

i) $\Pi$ is a leaf-tangent Poisson structure on $(M,\mathcal{F})$,
i.e., its restriction to each leaf is a Poisson structure of the
leaf;

ii) $d'\sigma=0$, equivalently, $d\sigma(X_1,X_2,X_3)=0$, $\forall
X_{1},X_{2},X_{3}\in\Gamma H$;

iii)  for any projectable (to the space of leaves of $\mathcal{F}$)
vector fields $X_1,X_2\in\Gamma_{pr}H$ ($pr$ denotes projectability)
one has
$$pr_F[X_1,X_2]=\sharp_{\Pi}(d''(\sigma(X_1,X_2)));$$

iv) for any projectable vector field $X\in\Gamma_{pr}H$ one has $
L_X\Pi=0$.
\end{prop} \begin{proof}
Condition i) is the equivalent of the fourth formula
(\ref{integrabacD}) if $L_F=L(\Pi)$ is the graph of $\Pi$. In the
coupling case, if we put $\alpha_a=\flat_\sigma X_a$ $(a=1,2,3)$ in
the first formula (\ref{integrabacD}), we get condition ii). The
similar replacement of the forms $\alpha$ in the second formula
(\ref{integrabacD}) puts the latter into the form
\begin{equation}\label{iiidincazDc}
(L_Y\sigma)(X_1,X_2)=-\lambda([X_1,X_2]),\; \forall
Y=\sharp_{\Pi}\lambda,\forall\lambda\in F^*,\forall X_1,X_2\in\Gamma
H.
\end{equation} Since this condition is invariant by
multiplication of the arguments $X$ by any $f\in C^\infty(M)$, it
suffices to ask (\ref{iiidincazDc}) for projectable arguments. But,
$X\in\Gamma_{pr}H$ iff $[Y,X]\in\Gamma F$, $\forall Y\in\Gamma F$
\cite{Mol}, and we see that (\ref{iiidincazDc}) is equivalent with
$$(\sharp_{\Pi}\lambda)(\sigma(X_1,X_2))=-\lambda([X_1,X_2]),$$
which exactly is condition iii) of the proposition. Similarly, it
suffices to use a projectable argument $X$ in the third formula
(\ref{integrabacD}). Then, the third formula (\ref{integrabacD})
becomes $([X,\sharp_{\Pi}\lambda],L_X\lambda)\in\Gamma L(\Pi)$,
which is equivalent with condition iv).
\end{proof} \begin{rem} \label{obsptDW} {\rm Conditions i)-iv) of Proposition
\ref{integrabcugd} are the same as Vorobiev's conditions
\cite{{Vor},{V04'}} of the Poisson case, except for the fact that
the $2$-form $\sigma$ may degenerate. In \cite{DW} these conditions
were included by definition. In the presymplectic case, where the
structure is defined by the closed $2$-form $\tau$ of
(\ref{descompsigma}) with a non degenerate component $\tau''$,
condition i) means that $\tau''$ defines symplectic structures on
the leaves of $\mathcal{F}$, and conditions ii) - iv) become
\begin{equation}\label{integrabpersympl} d'\tau'=0,\;
d''(\tau'(X_1,X_2))=\flat_{\tau''}(pr_F[X_1,X_2]),\;
L_X\tau''=0,\end{equation} $\forall X,X_1,X_2\in\Gamma_{pr}H$. It is
easy to see that these conditions are equivalent with
\begin{equation}\label{integrabpersympl2}d''\tau''=0,\;d'\tau'=0,
\;d''\tau'+\partial\tau''=0,\;d'\tau''=0,\end{equation} which are
the homogeneous components of $d\tau=0$.}
\end{rem}
\section{Dirac structures near a presymplectic leaf} We continue to
use the notation of the previous sections. A presymplectic leaf of a
Dirac structure $L$ of a differentiable manifold $M$ is an integral
submanifold of the distribution $L_{+}$ defined in Section 1. In
\cite{DW}, it was proven that, in a tubular neighborhood of an
embedded presymplectic leaf,  any Dirac structure $L$ is coupling
with respect to the fibers of the tubular structure. This result,
which extends the similar one in Poisson geometry \cite{Vor},
describes the geometry of a Dirac structure near an embedded
presymplectic leaf. Below, we give an invariant proof of this
result.
\begin{prop}\label{foaietransversa} Let $L$ be a Dirac structure
on the foliated manifold $(M,\mathcal{F})$. Assume that $L$ has a
presymplectic leaf $S$ such that $T_SM=TS\oplus F|_S$. Then, there
exists an open neighborhood $U$ of $S$ in $M$ such that $L|_U$ is
coupling with respect to $\mathcal{F}\cap U$.
\end{prop}
\begin{proof} We will refer to bidegrees defined
by the decomposition $T_SM=TS\oplus F|_N$, where we know that
$TS=L_{+S}$. Then, for all $X\in TS$, there exists a covector
$\alpha\in T^*S$ of bidegree $(1,0)$ equal to $\flat_{\omega^L_+}X$
on $TS$, and (\ref{reconstrincazclasic}) shows that $(X,\alpha)\in
L$. The conclusion is that $L|_S$ is the presymplectic structure of
$S$ and $H(L,\mathcal{F})|_S=TS$, therefore, $L$ is
$\mathcal{F}$-coupling along $S$. By Corollary \ref{Fernandes}, this
is equivalent with
\begin{equation}\label{ccusuma}
[L+(F\oplus ann\,F)]_S=T_SM\oplus T_S^*M,\end{equation} and, since
its left hand side is differentiable, (\ref{ccusuma}) also holds on
an open neighborhood $U$ of $S$.\end{proof}
\begin{corol}\label{corolarDW} {\rm \cite{DW}} Assume that the
Dirac structure $L$ of a manifold $M$ has an embedded leaf $S$.
Then, on a sufficiently small tubular neighborhood $U$ of $S$, with
the foliation $\mathcal{F}$ defined by the fibers of the tubular
structure, $L$ may be put under the form {\rm(\ref{Lreconstruit})},
where $(H,\sigma\in\Gamma\wedge^2(ann\,F),\Pi\in\Gamma\wedge^2F)$ is
a triple of geometric data that satisfies the integrability
conditions i)-iv) of Proposition {\rm \ref{integrabcugd}}, and $\Pi$
is a Poisson structure which vanishes on $S$.\end{corol}
\begin{rem} {\rm\cite{DW} The previous corollary implies the
fact that all the presymplectic leaves of a Dirac structure have the
same parity. Indeed, with the notation of the corollary, at a point
near $S$, $L_{+}$ is a direct sum of a subspace of dimension
$dim\,S$ and a subspace tangent to a symplectic leaf of $\Pi$, which
has an even dimension.} \end{rem}

Following \cite{DW}, we say that $\Pi$ is the {\it transversal
Poisson structure} of the leaf $S$. The following proposition shows
that the transversal structure is essentially unique.
\begin{prop}\label{unictransvers} {\rm\cite{DW}}
The transverse Poisson structure of an embedded presymplectic leaf
of a Dirac structure is unique up to Poisson equivalence.\end{prop}
\begin{proof} We get two transversal structures $\Pi_1,\Pi_2$ if we
use two tubular neighborhoods $U_1,U_2$. The isotopy of the latter
\cite{Hirsch} yields a leaf-preserving diffeomorphism
$\Phi:U_1\rightarrow U_2$, which may be seen as the composition of
maps in the flows of projectable vector fields $X$ on $U_1$. But,
the  (\ref{cond3inDW}) implies $L_X\Pi=0$ for all the projectable
vector fields. Hence, by the integrability condition iv) of
Proposition \ref{integrabcugd}, $\Phi_{*}\Pi_1=\Pi_2$ (after
shrinking the tubular neighborhoods if necessary).
\end{proof}

In what follows we consider the Vorobiev-Poisson structure
\cite{{Vor},{V04'}} in the case of an embedded, presymplectic leaf
$S$ of a Dirac structure $L$.

From the general result on Lie algebroids given by Theorem 2.1 of
\cite{IKV}, it follows that the vector bundle $L|_S$ has a well
defined induced structure of a transitive Lie algebroid over $S$
with the anchor and bracket defined by
\begin{equation}\label{algebroid}
\begin{array}{c} \rho(X,\theta)=X,\vspace{2mm}\\

[(X_1,\theta_1),(X_2,\theta_2)]_S=([\tilde X_1,\tilde X_2],L_{\tilde
X_1}\tilde\theta_2 -L_{\tilde
X_2}\tilde\theta_1-d(\tilde\theta_2(\tilde
X_1)))|_S,\end{array}\end{equation} where all the pairs $(X,\theta)$
belong to $L|_S$ and $(\tilde X,\tilde\theta)$ are arbitrary
extensions of $(X,\theta)$ from $S$ to $M$. The existence of such
extensions follows from the fact that $S$ is an embedded
submanifold, and the independence of the bracket (\ref{algebroid})
of the choice of the extensions follows easily if the extensions are
expressed via a local basis of $L$ and the axioms of a Lie algebroid
are used \cite{IKV}.

Accordingly, $\mathbb{G}=ker\,\rho$ is a bundle of Lie algebras such
that
\begin{equation}\label{sirexact1} 0\rightarrow\mathbb{G}
\stackrel{\subseteq}{\rightarrow}L|_S\stackrel{\rho}{\rightarrow}
TS\rightarrow0 \end{equation} is an exact sequence with projections
$p$ onto $S$, while $S$ is endowed with the presymplectic (closed)
$2$-form $\varpi=\omega_{+}^L$.

As in \cite{{Vor},{V04'}}, each splitting $\gamma: TS \rightarrow
L|_S$ produces geometric data $(\mathcal{H},\sigma,\mathbb{L})$ on
the total space of the dual bundle $\mathbb{G}^*$ with the foliation
by fibers $\mathcal{V}$. Namely: i) $\mathcal{H}$ is the horizontal
bundle of the dual of the connection defined on $\mathbb{G}$ by the
formula
\begin{equation}\label{conexpeG}
\nabla_X\eta=[\gamma(X),\eta]_S,\hspace{5mm}X\in\Gamma
TS,\,\eta\in\Gamma\mathbb{G},\end{equation} ii) $\sigma$ is the
$2$-form evaluated, at $z\in\mathbb{G}^*$, on horizontal lifts
$\mathcal{X}_1,\mathcal{X}_2$ of $X_1,X_2\in\Gamma TS$, by
\begin{equation}\label{sigmaVorob}
\sigma_z(\mathcal{X}_1,\mathcal{X}_2)=\varpi_{p(z)}(X_1,X_2) -
z([\gamma(X_1),\gamma(X_2)]-\gamma[X_1,X_2]),\end{equation} iii)
$\mathbb{L}$ is the family of Lie-Poisson structures of the fibers
of $\mathbb{G}^*$.

The only difference between this situation and that of Vorobiev
\cite{Vor} is that $\varpi$ may degenerate. But, the arguments and
computations of Vorobiev's case, as described in \cite{V04'} are
still valid, and they show that the triple
$(\mathcal{H},\sigma,\mathbb{L})$ satisfies the integrability
conditions i)-iv) of Proposition \ref{integrabcugd}. Therefore,
there exists a corresponding coupling Dirac structure
$\mathcal{L}(S,\gamma)$ on the manifold $\mathbb{G}^*$, and we call
it an {\it associated Dirac structure of $L$ along $S$}.

Let $\nu S$ be a normal bundle of the leaf $S$ $(T_SM=\nu S\oplus
TS)$. Then, the reconstruction formula (\ref{reconstrincazclasic})
shows that
\begin{equation}\label{LScuNS} L|_S=\{(X,\flat_\varpi
X+\lambda)\,/\,X\in TS,\,\lambda\in \nu^*S\},\end{equation} where
$\flat_\varpi X$ is the $1$-form defined by
\begin{equation}\label{flatptNS} (\flat_\varpi X)(Z)=
\left\{ \begin{array}{cl} \varpi(X,Z),&{\rm if}Z\in
TS,\vspace{2mm}\\ 0,&{\rm if}Z\in\nu S,\end{array}\right.
\end{equation} and we have a splitting $\gamma$ given
by \begin{equation}\label{splitting} \gamma(X)=(X,\flat_\varpi X).
\end{equation}

On the other hand, (\ref{LScuNS}) shows that we may identify the
bundle $\mathbb{G}^*$ with $\nu S$. Namely, we have
$\mathbb{G}=\nu^*S=ann\,TS$, and $\Theta\in Hom(ann\,TS,\mathbb{R})$
identifies with $Y\in \nu S$ defined by
$\lambda(Y)=\Theta(\lambda)$, $\forall\lambda\in ann\,TS.$ We will
say that the Dirac structure associated with $L$ by the splitting
(\ref{splitting}) is the {\it associated, normal Dirac structure}
$\mathcal{L}(S,\nu S)$.

We want to find a convenient local coordinate expression of
$\mathcal{L}(S,\nu S)$. For this purpose, around the points of $S$
and after the choice of the normal bundle $\nu S$, we take local
coordinates $(x^u,y^a)$, where $a$ (and similar indices $b,c$) takes
the values $1,...,codim\,S$ and $u$ (and similar indices $v,w$)
takes the values $1,...,dim\,S$, such that $S$ is locally defined by
the equations $y^a=0$ and
\begin{equation}\label{generatoripeS} \begin{array}{ll}
TS=span\left\{\frac{\partial}{\partial
x^u}|_{y^b=0}\right\},&\hspace*{1cm} \nu
S=span\left\{\frac{\partial}{\partial
y^a}|_{y^b=0}\right\},\vspace{2mm}\\ T^*S=span\{dx^u|_{y^b=0}\},&
\hspace*{1cm}\nu^*S=span\{dy^a|_{y^b=0}\}.\end{array} \end{equation}

Then, Theorem 3.2 of \cite{DW} tells us that the Dirac structure $L$
has local bases of the form
\begin{equation}\label{bazeleDW} \begin{array}{l}
\mathcal{H}_u=\left(\frac{\partial}{\partial x^u} + A_u^b(x,y)
\frac{\partial}{\partial
y^b},\alpha_{uv}(x,y)dx^v\right),\vspace{2mm}\\
\mathcal{V}^a=\left(B^{ab}(x,y) \frac{\partial}{\partial
y^b},dy^a-A^a_v(x,y)dx^v\right),\end{array}\end{equation} where
$\alpha_{uv}(x,0)$ are the components of the $2$-form $\varpi$ and
\begin{equation}\label{anularepeS}A^b_u(x,0)=0,\;B^{ab}(x,0)=0.
\end{equation} \begin{rem}\label{DWalgebric} {\rm As a matter of
fact the proof given in \cite{DW} holds to show that, for any
maximal isotropic subspace $L\subseteq W\oplus W^*$ where $W$ is an
arbitrary vector space, there are bases of the form
\begin{equation}\label{bazeDWalgebrice}
(l_u+A^b_uf_b,\alpha_{uv}\lambda^v),\,(B^{ab}f_b,\varphi^a-A^a_v\lambda^v),
\end{equation} where $(l_u)$ is a basis of $L_{+}$,
$(f_a)$ is a basis of an arbitrary complement of $L_{+}$,
$(\lambda^u)$ is the dual basis of $(l_u)$ and $(\varphi^a)$ is the
dual basis of $(f_a)$.}\end{rem}

It follows that the local basis of the Lie algebra bundle
$\mathbb{G}\approx \nu^*S$ given in (\ref{generatoripeS}) may be
seen as $(\mathcal{V}^a|_{y^b=0})$ and the local basis of
$\gamma(TS)$ is $(\mathcal{H}_u|_{y^b=0})$. Using
(\ref{anularepeS}), which implies the vanishing of the derivatives
of the same functions with respect to $x^u$ on $S$, and the
closedness of the $2$-form $\varpi$, we get for the brackets of the
elements of these bases the expressions
\begin{equation}\label{crosetebazice}\begin{array}{c}
[\mathcal{V}^a,\mathcal{V}^b]_{S}=\left(\frac{\partial
B^{ab}}{\partial y^c}\mathcal{V}^c\right)_{y^b=0},\,
[\mathcal{H}_u,\mathcal{V}^a]_{S}=\left(\frac{\partial
A^a_u}{\partial y^c}\mathcal{V}^c\right)_{y^b=0},\vspace{2mm}\\

[\mathcal{H}_u,\mathcal{H}_v]_{S}=\left(\frac{\partial
\alpha_{uv}}{\partial
y^c}\mathcal{V}^c\right)_{y^b=0}.\end{array}\end{equation}

Accordingly, as in \cite{V04'}, we get the following expressions of
the geometric data that define the Dirac structure
$\mathcal{L}(S,\nu S)$
\begin{equation}\label{orizontalptnormalL} \mathcal{H}=
span\{\mathcal{X}_u=\frac{\partial}{\partial x^u} + \frac{\partial
A^a_u}{\partial
y^c}\eta^c\frac{\partial}{\partial\eta^a}\}_{y^b=0},\end{equation}
\begin{equation}\label{sigmaptnormalL}
\sigma(\mathcal{X}_u,\mathcal{X}_v)=\alpha_{uv}(x,0)-
\left(\frac{\partial\alpha_{uv}}{\partial
y^a}\right)_{y^b=0}\eta^a,\end{equation}
\begin{equation}\label{bbLptnormalL}
\mathbb{L}^{ab}=\left(\frac{\partial B^{ab}}{\partial
y^c}\right)_{y^b=0} \eta^c,\end{equation} where $(\eta^a)$ are fiber
coordinates on $\nu S$.

If we use the previous formulas and (\ref{Lreconstruit}) we get a
canonical, local basis of type (\ref{bazeleDW}) for
$\mathcal{L}(S,\nu S)$ namely,
\begin{equation}\label{bazaptnormal}\begin{array}{l}
\left([\frac{\partial}{\partial x^u} + \frac{\partial
A^a_u}{\partial y^c}\eta^c\frac{\partial}{\partial\eta^a}]_{y^b=0},
[(\alpha_{uv}(x,0)- \left(\frac{\partial\alpha_{uv}}{\partial
y^a}\right)\eta^a)dx^v]_{y^b=0}\right),\vspace{2mm}\\

\left([\frac{\partial B^{ab}}{\partial y^c}\eta^c
\frac{\partial}{\partial\eta^b}]_{y^b=0}, [d\eta^a- \frac{\partial
A^a_u}{\partial y^c}\eta^cdx^u]_{y^b=0}\right).
\end{array}\end{equation} From (\ref{bazaptnormal}), we see that
$S$, seen as the zero section of $\nu S$ endowed with the $2$-form
$\varpi$, is a presymplectic leaf of $\mathcal{L}(S,\nu S)$.
Moreover, we see that the structure $\mathcal{L}(S,\nu S)$ along $S$
is a {\it linear approximation} of the Dirac structure $L$ along
$S$.

For a Poisson structure $P$, Vorobiev proved that the Poisson
structures $\mathcal{L}(S,\nu S)$ defined by different normal
bundles $\nu S$ are equivalent in neighborhoods of $S$
\cite{{Vor},{V04'}}. We will see below why his proof does not apply
in the general Dirac case, and indicate a particular case where it
works.

The choice of the normal bundle $\nu S$ is equivalent with the
definition of a projection epimorphism $\pi:T_SM\rightarrow TS$
$(\pi^2=\pi)$ namely, $\nu S=ker\,\pi$. A second normal bundle
$\nu'S$ corresponds to a second epimorphism $\pi':T_SM\rightarrow
TS$, and $\pi_t=(1-t)\pi +t\pi'$, $t\in\mathbb{R}$, defines a
homotopy between the normal bundles $\nu S$ $(t=0)$ and $\nu'S$
$(t=1)$ by a family of normal bundles $\nu_tS$.

Furthermore, there exists a bundle isomorphism $\Phi_1:\nu
S\rightarrow \nu'S$ defined by
$$\Phi_1(Y)=Y- \pi'(Y),\hspace{1cm}Y\in \nu S.$$ Similarly, we
have isomorphisms $$\Phi_t=t\Phi_1-(1-t)Id:\nu S\rightarrow
\nu_tS.$$ On each $\nu_tS$ we have a Dirac structure
$\mathcal{L}_t(S,\nu_tS)$, and we may pullback all these structures
to $\nu S$ by $\Phi_t^{-1}$.

On $\mathbb{G}^*$ seen as the fixed normal bundle $\nu S$, this
exactly provides the homotopy considered by Vorobiev
\cite{{Vor},{V04'}} between the Dirac structures
$\mathcal{L}(S,\gamma)$, $(\Phi^{-1}_1)_*(\mathcal{L}(S,\gamma'))$,
where $\gamma,\gamma'$ are the splittings of the exact sequence
(\ref{sirexact1}) associated with $\nu S,\nu'S$ by
(\ref{splitting}). Indeed, formula (\ref{splitting}) implies that
the $\mathbb{G}$-valued difference form $\phi=\gamma'-\gamma$ is
given by
\begin{equation}\label{formaphi}
\phi(X)=(i(X)\varpi)\circ\pi'\in ann\,TS\hspace{5mm}(X\in TS),
\end{equation} and, if we write the same form for $\pi_t$ instead of
$\pi'$ we get the form $t\phi$.

In the Poisson case, Vorobiev's proof is based on the $1$-form
\begin{equation}\label{formapsi}
\psi_Y(\mathcal{X})=<\phi(X),Y>=\varpi(X,\pi'Y),\hspace{5mm}Y\in \nu
S=\mathbb{G}^*,\,X\in TS,\end{equation} defined on the total space
of the bundle $\nu S$, where $\mathcal{X}$ is the horizontal lift of
$X$ by the connection (\ref{conexpeG}). The horizontal, time
dependent, tangent vector field $\Xi_t$ of $\nu S$ that satisfies
the condition $\flat_{\sigma_t}\Xi_t=-\psi$ has a flow which, at
time $1$, yields the required equivalence \cite{{Vor},{V04'}}.

In the general Dirac case, the form $\psi$ exists but, the vector
field $\Xi_t$ may not exist since the form $\sigma_t$ is no more non
degenerate.

Let us refer to the particular case of a Dirac structure $L$ such
that the field of planes $K=L\cap TM$ has a constant dimension $k$.
Since $L$ is closed by the Courant bracket, $K$ is involutive,
therefore, tangent to a foliation $\mathcal{K}$, the leaves of which
are submanifolds of the presymplectic leaves of $L$. A Dirac
structure with the previous property will be called {\it locally
reducible}, because if the stronger property that the leaves of
$\mathcal {K}$ are the fibers of a fibration holds the Dirac
structure is {\it reducible} \cite{LWX}. For a locally reducible
Dirac structure the local coordinates $(x^u)$ used to get the
canonical bases of $L$ given by (\ref{bazeleDW}) may be taken under
the form $(x^e,z^s)$ $(s=1,..,k,\,e=1,...dim\,S-k)$, where
$x^e=const.$ define the leaves of $\mathcal{K}$ and $(z^s)$ are
coordinates along these leaves.
\begin{prop}\label{bazeptstrred} Let $L$ be a locally reducible
Dirac structure, $x\in M$ and $S$ the presymplectic leaf through
$x$. Then, there exists a neighborhood of $x$ where $L$ has local
bases of the form
\begin{equation}\label{bazeleDWred} \begin{array}{l}
\mathcal{Z}_s=\left(\frac{\partial}{\partial z^s},0\right), \vspace{2mm}\\
\mathcal{H}_e=\left(\frac{\partial}{\partial x^e} + A_e^b(x,y)
\frac{\partial}{\partial
y^b},\alpha_{ef}(x,y)dx^f\right),\vspace{2mm}\\
\mathcal{V}^a=\left(B^{ab}(x,y) \frac{\partial}{\partial
y^b},dy^a-A^a_f(x,y)dx^f\right),\end{array}\end{equation} the
coefficients $A,B$ vanish at $y^a=0$ and $\alpha_{ef}(x,0)$ are the
local components of the presymplectic $2$-form $\varpi$ of $S$.
\end{prop}
\begin{proof} By the definition of $K$ and of the local
coordinates that we use, $\mathcal{Z}_s\in L$ and must be
$g$-orthogonal to the vectors of the basis (\ref{bazeleDW}). This
happens iff (\ref{bazeleDW}) are of the form (\ref{bazeleDWred}),
where, a priori, the remaining coefficients may depend on all the
coordinates $(x,y,z)$. But, since $L$ is closed by Courant brackets,
$[\mathcal{Z}_s,\mathcal{H}_e],[\mathcal{Z}_s,\mathcal{V}^a]$ must
be $g$-orthogonal to $\mathcal{Z}_s,\mathcal{H}_e$ and this implies
$$\frac{\partial\alpha_{ef}}{\partial z^s}=0,\,\frac{\partial
A_e^b}{\partial z^s}=0,\,\frac{\partial B^{ab}}{\partial z^s}=0.$$
\end{proof}

Now, it is easy to prove the following proposition:
\begin{prop}\label{invarincazred} Let $L$ be a reducible
Dirac structure and $S$ an embedded presymplectic leaf of $L$. Then,
the associated, normal Dirac structures defined by any two normal
bundles of $S$ are equivalent.
\end{prop}
\begin{proof} If $\psi:M\rightarrow M/\mathcal{K}$ is the
reducibility fibration, the vectors $\mathcal{H}_e,\mathcal{V}^a$ of
the canonical bases (\ref{bazeleDWred}) are $\psi$-projectable and
their projections define a Poisson structure $\Lambda$ on
$M/\mathcal{K}$ for which $S/\mathcal{K}$ is an embedded symplectic
leaf. (This is a well known result \cite{{C},{LWX}}. The structure
$\Lambda$ is Poisson because, on $S$, $K=ker\,\varpi$, hence, the
matrix $(\alpha_{ef})$ of (\ref{bazeleDWred}) has a maximal rank.)
Moreover, the projections of the vectors $\mathcal{V}^a$ of
(\ref{bazeleDWred}) yield a normal bundle of $S/\mathcal{K}\subseteq
M/\mathcal{K}$ and, because the computation of the Courant brackets
of $\mathcal{H}_e,\mathcal{V}^a$ on $M$ and on $M/\mathcal{K}$ is
the same, the associated, normal Dirac structure of $L$ along $S$
and that of the projected Poisson structure $\Lambda$ of
$M/\mathcal{K}$ along $S/\mathcal{K}$ correspond each to the other
by $\psi$. If we act like that for two normal bundles of $S$, we get
two associated, normal Poisson structures of $\Lambda$, which are
Poisson equivalent by Vorobiev's theorem \cite{{Vor},{V04'}}. This
Poisson equivalence lifts to an equivalence of the associated,
normal Dirac structures of $L$. (The triples that define the
associated structures on $M$ and on $M/\mathcal{K}$ have the same
local coordinate expressions with respect to the bases
(\ref{bazeleDWred}) and their projections.) \end{proof}
\section{Submanifolds of Dirac manifolds} In this section we will
show that the almost-coupling and coupling conditions are also
significant for single submanifolds $N^p$ of a Dirac manifold
$(M^n,L)$ (a differentiable manifold $M$ with a fixed Dirac
structure $L$), as opposed to a whole foliation $\mathcal{F}$. For
simplicity, we will assume that the submanifold is embedded. Where
possible, we will continue to use the notation of the previous
sections.

We begin by recalling that Dirac structures may be both pulled
back and pushed forward pointwisely (e.g., see \cite{{C},{BC}}).
If $\phi:N\rightarrow M$ is a differentiable mapping between
arbitrary manifolds and $L(M)$ is a Dirac structure on $M$ then,
$\forall x\in N$,
\begin{equation}\label{strDinapoi} (\phi^*(L(M)))_x = \{(Z,\phi^*\alpha)
\,/\,Z\in T_xM\,,\alpha\in T^*_{\phi(x)}M,\,\end{equation}
$$(\phi_{*,x}Z,\alpha)\in L(M)_{\phi(x)}\}$$
is a maximal isotropic subspace of $T_xN\oplus T^*_xN$

On the other hand, if we have a Dirac structure $L(N)$ on $N$ and
$x\in N$, we have the maximal isotropic subspace
\begin{equation}\label{strDinainte} (\phi_{*}(L(N)))_{\phi(x)} =
\{(\phi_{*,x}Z,\alpha)\,/\,\alpha\in
T^*_{\phi(x)}M,\,(Z,\phi^*_x\alpha)\in L(N)_x\}\end{equation}
$$\subseteq T_{\phi(x)}M\oplus T^*_{\phi(x)}M.$$

Generally, these pointwise operations do not yield differentiable
subbundles. If differentiable Dirac structures $L(N),L(M)$ are
related by (\ref{strDinapoi}), $\phi$ is called a {\it backward
Dirac map}, and if the relation is (\ref{strDinainte}) $\phi$ is a
{\it forward Dirac map} \cite{BC}. If $\phi$ is the embedding
$\iota:N^p\hookrightarrow M^n$ of a submanifold and if
$L(N)=\iota^*(L(M))$ is differentiable, $L(N)$ must be integrable,
and it defines a Dirac structure on $N$ \cite{C}. Indeed, $L(N)$ is
equivalent in the sense of (\ref{reconstrL+}) with the field of
planes
\begin{equation}\label{LNplus} L(N)_{+,x}=\{Z\in T_xN\,/\,
\exists\alpha\in T_x^*M,\,(Z,\alpha)\in L(M)_x\}\end{equation}
$$=L(M)_{+,x}\cap T_xN= T_x(S(L(M)))\cap T_xN$$ ($S$
denotes presymplectic leaves), endowed with the $2$-form
$\omega^{L(N)}_+$ induced by $\omega_{+}^{L(M)}$. Obviously, if
(\ref{LNplus}) is a differentiable distribution, it is integrable
and $\omega^{L(N)}_+$ is closed. In what follows, if the {\it
induced Dirac structure} $L(N)=\iota^*(L)$ is differentiable, we
will call $N$ a {\it proper submanifold} of $(M,L)$.

Along the submanifold $N$ of $(M,L)$, we have the field of planes
\begin{equation}\label{defK} K(N)=ker\,\flat_{\omega_{S(N)}}=L(N)\cap TN
\end{equation} $$=\{(Z,0)\,/\,Z\in TN\,\&\, \exists\alpha\in ann\,
TN,\,(Z,\alpha)\in L\}$$ $$=pr_{TN}(L\cap(TN\oplus ann\,TN)),$$
the {\it kernel} of the induced structure $L(N)$. If $L$ is
defined by a Poisson bivector field $P\in\mathcal{V}^2(M)$, a
proper submanifold with kernel zero has an induced structure
$L(N)$ provided by a Poisson bivector field
$\Pi\in\mathcal{V}^2(N)$. Such submanifolds were studied in
\cite{CF1} under the name of {\it Poisson-Dirac submanifolds}.

On the other hand, if $\iota:N\hookrightarrow M$ is a submanifold of
$(M,L)$ and $\nu N$ is a normal bundle, we may use the push forward
construction  (\ref{strDinainte}) along $N$ and get a maximally
isotropic subbundle $pr_{N}L\subseteq TN\oplus T^*N$ given by
\begin{equation}\label{proiectiastrD} pr_NL =
\{(pr_{TN}Z,pr_{T^*N}\alpha)\,/\,(Z,\alpha)\in L|_N\},
\end{equation} where the involved projections are those
of the decomposition $T_NM=\nu N\oplus TN$. Obviously, this
subbundle is differentiable. A pair $(N,\nu N)$ is called a {\it
normalized submanifold} \cite{V2}, and, along it, one has {\it
adapted local coordinates} of $M$ like those that appear in
(\ref{generatoripeS}). In the present case, if $pr_NL$ of
(\ref{proiectiastrD}) is also integrable, it yields a Dirac
structure and we will say that $(N,\nu N,pr_NL)$ is a {\it submersed
submanifold} of $(M,L)$.

Now, like in Definition \ref{precD}, we define
\begin{defin}\label{properlyn} {\rm The pair $(N,\nu N)$ is a {\it
properly normalized submanifold} of the Dirac manifold $(M,L)$ if
\begin{equation}\label{descpropnorm} L|_N=L_{\nu N}\oplus L_{TN},
\end{equation} where} \begin{equation}\label{compoluiLindescn}
L_{\nu N}=L\cap(\nu N\oplus\nu^*N),\,L_{TN}=L\cap(TN\oplus T^*N).
\end{equation}\end{defin}
\begin{prop} A properly normalized submanifold $(N,\nu N)$ of a
Dirac manifold $(M,L)$ is simultaneously proper and submersed and
has the differentiable, Dirac structure
\begin{equation} \label{strdinproperlyn}
L(N)=\iota^*(L)=pr_{N}L,\end{equation} where the projection is
defined by the decomposition $T_NM=\nu N\oplus TN$. \end{prop}
\begin{proof} It is easy to understand that condition (\ref{descpropnorm}) is
equivalent with
\begin{equation}\label{eqproperlyn} (Z,\theta)\in L|_N \Rightarrow
(pr_{TN}Z,pr_{T^*N}\theta)\in L|_N. \end{equation} Accordingly,
the pair defined by $(Z,\theta)\in L$ in $\iota^*(L)$ is the same
as that defined by $(Z,pr_{T^*N}\theta)$ and may be identified
with the latter. This justifies the equalities
(\ref{strdinproperlyn}) and proves the proposition. \end{proof}

\begin{rem}\label{Dprojection} {\rm If $L=L(P)$
where $P$ is a Poisson bivector field, condition
(\ref{descpropnorm}) reduces to the almost-coupling condition of
the Poisson case (see the beginning of Section 2 or \cite{V04'})
and the manifold is a {\it Poisson-Dirac submanifold with a Dirac
projection} in the sense of
\cite{CF1}.} \end{rem}

Furthermore, along any submanifold $N$ of $(M,L)$ we may define a
field of subspaces $H_x(L,N)$, $x\in N$ by formula
(\ref{distribHLF}) with condition $\alpha\in ann\,F_x$ changed to
$\alpha\in ann\,T_xN$. This leads to a notion that corresponds to
coupling namely,
\begin{defin}\label{cosimplecticD} {\rm The submanifold $N$ is a
{\it cosymplectic submanifold} of the Dirac manifold $(M,L)$ if,
$\forall x\in N$, $T_xM=H_x(L,N)\oplus T_xN$.}\end{defin} The
reason for this name is that if $L=L(P)$ for a Poisson bivector
field $P$ then $N$ is a cosymplectic submanifold in the sense of
\cite{W1}. Obviously, now, all the results stated in Proposition
\ref{cimplicaac} are true along $N$. In particular, $L|_N$ has
the expression (\ref{Lreconstruit}) (with $F^*$ replaced by
$T^*N$) and the induced structure $L(N)$ is Poisson and defined by
the bivector field $\Pi$ of (\ref{Lreconstruit}).\\

In what follows, we define a restriction of the Courant bracket of
$L$ to a submanifold $\iota:N\hookrightarrow(M,L)$. Like for
Poisson-Dirac submanifolds \cite{CF1}, we may define
\begin{equation}\label{AptL} A_N(M,L)=\{(X,\alpha)\in L|_N\,/\,X\in
TN\},\end{equation} which is important because, by
(\ref{strDinapoi}), we have \begin{equation}\label{LNcuA}
L(N)=\iota^*(L)=\{(X,\iota^*\alpha)\,/\,(X,\alpha)\in A_N(M,L)\}.
\end{equation} Even though $A_N(M,L)$ may not be a vector bundle, we will
consider the real, linear space $\Gamma A_N(M,L)$ of differentiable
cross sections of $A_N(M,L)$ (which may be zero). Using a partition
of unity that consists of a tubular neighborhood of $N$ and open
sets that do not intersect $N$, it follows easily that any cross
section $(X,\alpha)\in\Gamma A_N(M,L)$ admits extensions $(\tilde
X,\tilde\alpha)\in \Gamma L$. Accordingly, on $\Gamma A_N(M,L)$ we
may define a bracket
\begin{equation}\label{crosetpeA} [(X,\alpha),(Y,\beta)]_A=
([(\tilde X,\tilde\alpha),(\tilde Y,\tilde \beta)]_L)|_N.
\end{equation}
\begin{prop}\label{propcrosetpeA} The bracket {\rm(\ref{crosetpeA})} is
well defined and, together with the projection $\rho(X,\alpha)=X$,
yields a structure of Herz-Reinhart Lie algebra over
$(\mathbb{R},C^\infty(N))$ on $\Gamma A_N(M,L)$.\end{prop}
\begin{proof} In order to prove that the bracket (\ref{crosetpeA})
does not depend on the choice of the extensions it suffices to prove
that it vanishes if, say, $(Y,\beta)=(0,0)$ (see the proof of
Theorem 2.1 of \cite{IKV}). If
$$(\tilde Y,\tilde \beta)=\sum_{i=1}^n\tilde\lambda_i(\tilde
B_i,\tilde\theta_i),$$ where $(\tilde B_i,\tilde\theta_i)$ is a
local basis of $L$ and $\tilde\lambda_i|_N=0$, and since
$X\in\mathcal{V}^1(N)$, we get
$$([(\tilde X,\tilde\alpha),(\tilde Y,\tilde \beta)]_L)|_N =
\sum_{i=1}^n(\tilde\lambda_i[(\tilde X,\tilde\alpha),(\tilde
B_i,\tilde\theta_i)]_L)|_N + \sum_{i=1}^n(X\tilde\lambda_i)(\tilde
B_i,\tilde\theta_i)|_N=0. $$

The last assertion of the proposition is obvious if we recall that a
Herz-Reinhart Lie (HRL) algebra (a pseudo-Lie algebra in the sense
of \cite{MacK}) $\mathcal{R}$ over $(\mathbb{R},C^\infty(N))$ is a
real Lie algebra, which is a $C^\infty(N)$-module endowed with a
homomorphism $\rho:\mathcal{R}\rightarrow \mathcal{V}^1(N)$, such
that the properties of the algebra of global cross sections of a Lie
algebroid hold.\end{proof}

Furthermore, the mapping $(X,\alpha)\mapsto(X,\iota^*\alpha)$
defines a HRL-algebra morphism
\begin{equation}\label{morfismHRL}
\iota^{\#}:\Gamma A_N(M,L)\rightarrow\Gamma(L(N))
\hspace{5mm}(L(N)=\iota^*(L))\end{equation}
and $ker\,\iota^{\#}=\Gamma(L\cap ann\,TN)$.

\begin{prop}\label{diffA} If $(N,\nu N)$ is a properly normalized
submanifold of $(M,L)$, $A_N(M,L)$ is a differentiable field of
subspaces of the fibers of $TN\oplus T^*_NM$. \end{prop}
\begin{proof} For any $(X,\alpha)\in A_N(M,L)$, we may write
\begin{equation} \label{descompaux1}
(X,\alpha)=\sum_{i=1}^n\lambda_i(pr_{TN}B_i,pr_{T^*N}\theta_i)
+(0,pr_{\nu^*N}\alpha),\end{equation} where $(B_i,\theta_i)$ is a
local basis of $L|_N$. Hence, the local cross sections of
$A_N(M,L)$ are spanned by differentiable cross sections.
\end{proof}

As a consequence of Proposition \ref{diffA} and of formula
(\ref{LNcuA}), in the case of a properly normalized submanifold
the morphism $\iota^{\#}$ is surjective, and we have the following
exact sequence of HRL-algebras
\begin{equation}\label{sirexactHRL} 0\rightarrow\Gamma(L\cap
ann\,TN)\rightarrow\Gamma
A_N(M,L)\rightarrow\Gamma(L(N))\rightarrow0.
\end{equation}\\

Along a properly normalized submanifold $(N,\nu N)$, all the
vector fields and differential forms split into $TN$ and $\nu
N$-components, and we may identify $\Gamma T^*N$ with the space of
the tangent components of $1$-forms in $\Gamma T^*_NM$.
Accordingly, we may write a decomposition formula
\begin{equation}\label{Gauss} [(X,\lambda),(Y,\mu)]_A
= [(X,\lambda),(Y,\mu)]_{L(N)}
+(0,B((X,\lambda),(Y,\mu))),\end{equation} where $\lambda,\mu\in
pr_{T*N}(T^*_NM)$, $B((X,\lambda),(Y,\mu))\in\nu^*N$ and $\forall
Z\in\nu N$ one has
\begin{equation}\label{defformadoi}
B((X,\lambda),(Y,\mu))(Z)= Z(\tilde\lambda(\tilde Y)) -
\lambda([\tilde Z,\tilde Y])+\mu([\tilde Z,\tilde X]),
\end{equation} where $\tilde Z$ is an extension of $Z$ to $M$. The
result does not depend on the choice of the extension $\tilde Z$
because $B((X,\lambda),(Y,\mu))$ is a $1$-form. By an analogy with
Riemannian geometry to be explained below, we call $B$ the {\it
second fundamental form} of $(N,\nu N)$.

Let $N$ be a Poisson-Dirac submanifold of the Poisson manifold
$(M,P)$, which is properly normalized by the normal bundle $\nu N$
and has the induced Poisson structure $\Pi$. Then, the kernel
condition $K(N)=0$ (see (\ref{defK})) becomes
\begin{equation}\label{cond1dePD} TN\cap\sharp_P(ann\,TN)=0,
\end{equation}or, by passing to the annihilator spaces,
\begin{equation}\label{cond2dePD} ann\,TN + A_N(M,P) =
T^*_NM,\end{equation} where
\begin{equation}\label{defluiA}A_N(M,P)=
\{\xi\in T^*_NM\,/\,\sharp_P\xi\in TN\}\end{equation}
$$ = ann(\sharp_P(ann\,TN))\approx A_N(M,L(P)).$$

Furthermore, the bracket (\ref{crosetpeA}) produces a bracket of
$1$-forms
\begin{equation}\label{crosetformeNcazPD}
\{\alpha,\beta\}_A=\{\tilde\alpha,\tilde\beta\}_P|_N\in\Gamma
T^*_NM\hspace{5mm}(\alpha,\beta\in \Gamma T^*N),\end{equation}
where the $P$-bracket is that of the cotangent Lie algebroid of
$(M,P)$ and tilde denotes extension to $M$. Formula (\ref{Gauss})
becomes
\begin{equation}\label{Gauss1}
\{\alpha,\beta\}_A=\{\alpha,\beta\}_\Pi +
B(\alpha,\beta),\end{equation} where the second fundamental form
$B$ is given by
\begin{equation}\label{formadouaptP} B(\alpha,\beta)(Z)= -(L_{\tilde Z}P)|_N
(\alpha,\beta),\hspace{5mm}Z=\tilde Z|_N\in\Gamma\nu
N.\end{equation}

Now, the Riemannian terminology used above is justified as follows.
A Riemannian metric $g$ of the Poisson manifold $(M,P)$ yields a
canonical {\it Riemannian, contravariant derivative} $D^P$ \cite{B}.
This is a cotangent-Lie algebroid-connection on $TM$ which preserves
the metric and has no torsion i.e.,
\begin{equation}\label{1conexR}
(\sharp_P\gamma)(g(\alpha,\beta))=g(D^P_\gamma\alpha,\beta)+g(\alpha,
D^P_\gamma\beta),\end{equation}
\begin{equation}\label{2conexR}
D^P_\alpha\beta-D^P_\beta\alpha=\{\alpha,\beta\}_P,\end{equation}
for all $\alpha,\beta,\gamma\in\Omega^1M$. This operator is provided
by the usual algebraic trick that derives the Riemannian connection
from the metric, and the result is
\begin{equation}\label{3conexR} \begin{array}{r}2g(D^P_\alpha\beta,\gamma)=
(\sharp_P\alpha)(g(\beta,\gamma))
+(\sharp_P\beta)(g(\gamma,\alpha))
-(\sharp_P\gamma)(g(\alpha,\beta))\vspace{2mm}\\ +
g(\{\alpha,\beta\}_P,\gamma) +g(\{\gamma,\alpha\}_P,\beta)+
g(\{\gamma,\beta\}_P,\alpha).\end{array}\end{equation}

Now, assume that $(N,\nu N=T^{\perp_g}N)$ is a properly normalized
Poisson-Dirac submanifold with the induced Poisson structure
$\Pi$. Then $N$ has its own canonical operator $D^\Pi$ on $T^*N$
and, also, a contravariant derivative $D^{P,N}$ defined by
\begin{equation}\label{conexPN} D^{P,N}_\alpha\beta =
(D^P_{\tilde\alpha}\tilde\beta)|_N,\hspace{5mm}
\alpha,\beta\in\Gamma T^*N,\end{equation} where
$\tilde\alpha,\tilde\beta$ are extensions of $\alpha,\beta$. It
follows easily from (\ref{3conexR}) that the result of
(\ref{conexPN}) does not depend on the choice of the extension.
The formula
\begin{equation}\label{Gauss2} D^{P,N}_\alpha\beta = D^{\Pi}_\alpha\beta
+
\Psi(\alpha,\beta),\hspace{5mm}Psi\in\Gamma\otimes^2TM,\end{equation}
is a {\it Gauss-type equation} and $\Psi$ is the $g$-{\it second
fundamental form} of $N$. But, (\ref{1conexR}) shows that $\Psi$ is
determined by the form $B$ of (\ref{Gauss1}). Namely, we get
\begin{equation}\label{ScuB} -2g(\Psi(\alpha,\beta),\gamma)=
g(B(\alpha,\beta),\gamma)+ g(B(\gamma,\alpha),\beta) +
g(B(\gamma,\beta),\alpha).\end{equation} The tensor field $\Psi$
is not symmetric and its skew-symmetric part is $(1/2)B$.
\begin{rem}\label{nupring} {\rm For any normalized submanifold
$(N,\nu N)$ of any manifold $M$, there exist Riemannian metrics $g$
of $M$ such that $\nu N=T^{\perp_g}N$. To get one, it suffices to
define it along $N$, then extend to $M$ along an open covering that
consists of a tubular neighborhood of $N$ and of sets that do not
intersect $N$ via a partition of unity. Then, if $M$ is endowed with
a Poisson structure $P$, it follows easily that $(N,\nu N)$ is a
properly normalized Poisson-Dirac submanifold iff $N$ is invariant
by $\Phi=\sharp_P\circ\flat_g$. In particular, if $M$ is a K\"ahler
manifold, $\Phi$ is the complex structure tensor and the properly
normalized Poisson-Dirac submanifolds are the complex analytic
submanifolds of $M$.}\end{rem}
\begin{prop}\label{Bptcoymplectic} A cosymplectic submanifold $N$ of
a Dirac manifold $(M,L)$ has a vanishing second fundamental form.
\end{prop} \begin{proof} For a cosymplectic submanifold $N$,
Corollary \ref{Fernandes} holds for $TN$ instead of $F$ and, in
particular, $L\cap ann\,TN=0$. Then, by (\ref{eqproperlyn}),
$(X,\alpha)\in A_N(M,L)$  implies $(0,pr_{\nu^*N}\alpha)\in L$,
hence, we must have $pr_{\nu^*N}\alpha=0$. Therefore,
$A_N(M,L)=L(N)$, $\iota^{\#}$ is an isomorphism, and $B=0$.
\end{proof}
\begin{defin}\label{deftotalD} {\rm A submanifold $N$ of a Dirac
manifold $(M,L)$ which has a normal bundle $\nu N$ such that $(N,\nu
N)$ is properly normalized and has a vanishing second fundamental
form will be called a {\it totally Dirac submanifold}.}\end{defin}

In the Poisson case, these submanifolds were called Dirac \cite{Xu1}
or Lie-Dirac \cite{CF1}. We took the term {\it totally} from
Riemannian geometry (totally geodesic submanifolds). \vspace*{1mm}\\

{\it Acknowledgements}. Part of the work on this paper was done
during a visit to the Centre Bernoulli of the \'Ecole Polytechnique
F\'ederale de Lausanne, Switzerland, the support of which is
acknowledged here, with thanks. The author is also grateful to Rui
Fernandes for a private communication that led to an improvement of
a previous version of the paper.
\hspace*{7.5cm}{\small \begin{tabular}{l} Department of
Mathematics\\ University of Haifa, Israel\\ E-mail:
vaisman@math.haifa.ac.il \end{tabular}} \hspace*{7.5cm}{\small
\begin{tabular}{l} Department of Mathematics\\ University of Haifa,
Israel\\ E-mail: vaisman@math.haifa.ac.il \end{tabular}}
\end{document}